\documentclass[12pt]{article}
\newcommand{\doublespace}{
   \renewcommand{\baselinestretch}{1.2}
   \large\normalsize}

\def \Z{\Bbb Z}
\def \C{\Bbb C} 
\def \R{\Bbb R}
\def \Q{\Bbb Q}
\def \N{\Bbb N}

\def \wt{{\rm wt}}

\def \span{{\rm span}}
\def \max{{\rm max}}
\def \wt{{\rm wt}}

\def \dim{{\rm dim}}
\def \Res{{\rm Res}}

\def \End{{\rm End}}

\def \<{\langle}
\def \>{\rangle}

\def \l{\lambda }

\def \om{\omega }

\def \pf{\noindent {\bf Proof:} \,}
\def \cg{\chi_g}
\def \cg'{\chi'_g}

\def \d{\delta}

\def \1{{\textbf 1}}

\def \qed{{\mbox{$\square$}}}
\input amssym.def
\input amssym
\doublespace
\begin{document}
\bibliographystyle{alpha}
\newtheorem{thm}{Theorem}[section]
\newtheorem{thmn}{Theorem}
\newtheorem{prop}[thm]{Proposition}
\newtheorem{cor}[thm]{Corollary}
\newtheorem{lem}[thm]{Lemma}
\newtheorem{rem}[thm]{Remark}
\newtheorem{de}[thm]{Definition}
\newtheorem{hy}[thm]{Hypothesis}
\begin{center}
{\Large {\bf Ordered spanning sets for vertex operator algebras and their modules}\footnote{A contribution to the Moonshine Conference at ICMS, Edinburgh, July 2004.}}
\\
\vspace{0.5cm} Geoffrey Buhl\footnote{Supported by a NSF Postdoctoral Fellowship for the Mathematical Sciences.}\\
Mathematics Department\\
       Rutgers University\\
       Piscataway, NJ 08854-8019\\
{\tt gbuhl@math.rutgers.edu}
\end{center}

\begin{abstract}
Moonshine relates three fundamental mathematical objects: the Monster sporadic simple group, the modular function $j(\tau)$, and the moonshine module vertex operator algebra $V^\natural$.  Examining the relationship between modular functions and the representation theory of vertex operator algebras reveals rich structure.  In particular, $C_2$-cofiniteness (also called Zhu's finiteness condition) implies the existence of finite generating sets and Poincar\'{e}-Birkhoff-Witt-like spanning sets for vertex operator algebras and their modules.  These spanning sets feature desirable ordering restrictions, e.g.,  a difference-one condition.
\end{abstract}

%%%%%%%%%%%%%%%%%%%%%%
\section{Introduction}
%%%%%%%%%%%%%%%%%%%%%%

The theory of vertex operator algebra blossomed from two major accomplishments: the proof of the McKay-Thompson conjecture by Frenkel, Lepowsky, and Meurman \cite{MR996026} who constructed the Moonshine module $V^\natural$ and the proof of  the Conway-Norton conjecture by Borcherds  \cite{MR1172696} using the Moonshine module.  These two conjectures make up what is commonly referred to as Monstrous Moonshine, relating the modular function $j(\tau)$ and the Monster group by way of a third fundamental mathematical object, the Moonshine module vertex operator algebra $V^\natural$.  The study of vertex operator algebras continues to reveal relations within mathematics and with physics. 

Representation theory is a particularly rich aspect of the theory of vertex operator algebras with fundamental connections to number theory, the theory of simple groups, and string and conformal field theories in physics.  A core idea in the representation theory of vertex operator algebras and conformal field theory is ``rationality'', a term used in a variety of ways to describe certain desirable properties of a vertex operator algebra and its modules. Complete reducibility of modules is one property that ``rationality'' invariably encompasses, but not always solely.  In both mathematics and physics, ``rationality'' is a term that suffers from a variety of meanings.  Compounding this difficulty is the variety of module definitions that appear in mathematics and physics literature.  This combination makes the concept ``rationality is complete reducibility of modules''  murky at best. 

For certain vertex operator algebras, we can achieve some clarity. An assumption on the ``size'' of the vertex operator algebra has important implications for its representation theory.  This size condition is called $C_2$-cofiniteness, and it implies the existence of a finite generating set and Poincar\'{e}-Birkhoff-Witt-like ordered spanning sets for the algebra and modules. An assumption of $C_2$-cofiniteness on an algebra ensures that even the most basic notion of a module has ``suitable'' structure.  In addition, the assumption of $C_2$-finiteness clarifies the concept of complete reducibility for modules of a vertex operator algebra.  Understanding the implications of $C_2$-cofiniteness is especially important in light of the recent developments in the representation theory of ``non-rational'', $C_2$-cofinite theories \cite{math.QA/0503472} \cite{carqueville-2006-39} \cite{math.QA/0311235}.

%%%%%%%%%%%%%%%%%%%%%%%
\section{Vertex operator algebras and quotient spaces}
%%%%%%%%%%%%%%%%%%%%%%%

For an introduction to the theory of vertex operator algebras, I refer the reader to ``Introduction to Vertex Operator Algebras and Their Representations'' by Lepowsky and Li \cite{MR2023933}.  Throughout this exposition, I will assume that the vertex operator algebras are of ``CFT-type''.  That is, a vertex operator algebra $V$ is of CFT-type if $V=\bigoplus_{n\geq 0} V_n$ and $V_0=\C\1$.  The weight of a homogenous vector is its $L(0)$-eigenvalue, $L(0)u= (\wt u) u$. The weight of an operator, or ``mode'', $u_n$ is also given by the $L(0)$-action, $L(0) u_n v = \wt(u_n)u_n L(0) v= (\wt u -n -1) u_n L(0) v$ for $n \in \Z$.

One of the powerful tools in the study of these infinite-dimesional objects, vertex operator algebras, has been to look at quotient spaces.  This technique's most important example is Zhu's algebra $A(V)$ \cite{MR1317233}.  

\begin{de} {\rm
For $V$ a vertex operator algebra, let $$O(V)=\span\{\Res_x \frac{(1+x)^{\wt u}}{x^{2}}Y(u,x)v:  u,v \in V \},$$ and let $A(V)=V/O(V)$.}
\end{de}

Zhu's algebra $A(V)$ is an associative algebra with identity, and it  acts on lowest weight vectors of modules.  This concept has been expanded to act on larger ``slices'' of modules.  The $n$th Zhu's algebra $A_n(V)$, also an associative algebra, acts on the bottom $n$ levels of modules \cite{MR1637252}. 

\begin{de} {\rm
For $V$ a vertex operator algebra and $n \in \N$, let $$O_n(V)=\span\{\Res_x \frac{(1+x)^{\wt u+n}}{x^{2n+2}}Y(u,x)v:  u,v \in V \},$$ and let $A_n(V)=V/O_n(V)$, where $A_0(V)=A(V)$.}
\end{de}

Under certain assumptions, all of the $n$th Zhu's algebras are semisimple and hence finite-dimensional.  The representation theories of a vertex operator algebra $V$ and its Zhu's algebras $A_n(V)$ are intimately related \cite{MR1317233} \cite{MR1615132} \cite{MR1637252}.

Another family of subspaces spaces used to create interesting quotient spaces are the $C_n$ spaces.  The subspace $C_2(V)=\span\{\Res_x x^{-2} Y(u,x)v: u,v \in V\}$ was introduced in Zhu's modularity paper \cite{MR1317233}.  One of the crucial assumptions needed to prove the modularity properties of certain graded traces is finite-dimensionality of the quotient space $V/C_2(V)$.  This property is known as $C_2$-cofiniteness or Zhu's finiteness condition.

This quotient space $V/C_2(V)$ has the structure of a Poisson algebra.  A Poisson algebra has two operations: an associative product $\cdot$ and a Lie bracket $[,]$ with compatibility of these operations given by Liebniz's Law $[x,y \cdot z]=[x,y]\cdot z+ y\cdot[x,z]$.  For $V/C_2(V)$, the product is given by $u \cdot v = \Res_x x^{-1}Y(u,x)v=u_{-1}v$ and the Poisson bracket is given by $[u,v]=\Res_x Y(u,x)v=u_0v$.  We can expand the definition of $C_2(V)$ to obtain a family of subspaces.  

\begin{de} {\rm
For a vertex operator algebra $V$ and for $n \geq 2$, let $$C_n(V)=\span\{\Res_x x^{-n} Y(u,x)v: u,v \in V\}.$$ Then $V$ is called {\it $C_n$-cofinite} if $V/C_n(V)$ is finite-dimensional.}
\end{de}

The case where $n=1$ is more nuanced and depending on an author's focus, is approached differently.  Focusing on the algebra, the naive extension of the definition, $\span\{\Res_x x^{-1} Y(u,x)v: u,v \in V\}$, is not particularly interesting since the creation axiom for vertex operator algebras ensures that this subspace is all of $V$.  A more interesting subspace is the following.

\begin{de} {\rm (cf. \cite{MR1700507})
For a vertex operator algebra $V=\bigoplus_{n \geq 0}V_n$, let $$C_1(V)=\span\{\Res_x x^{-1} Y(u,x)v, L(-1)u: u,v \in \bigoplus_{n>0}V_n\}.$$  Then $V$ is called {\it $C_1$-cofinite} if $V/C_1(V)$ is finite-dimensional.}
\end{de}

% How useful?
%For many of the quotient spaces described above, one can define analogous quotient spaces of modules.  

The assumption of $C_n$-cofiniteness of a vertex operator algebra controls the size of other quotient spaces.  For example, a simple calculation shows that if a vertex operator algebra $V$ is $C_2$-cofinite, then $A(V)$ is finite-dimensional. The $L(-1)$ derivation property implies $C_{n}(V) \subseteq C_{n-1}(V)$, and thus $C_n$-cofiniteness implies $C_{n-1}$-cofiniteness for $n \geq 2$.   In fact $C_2$-cofiniteness implies that a great deal of quotient spaces are finite-dimensional \cite{MR1990879}.

There are other interesting quotient spaces.  For example, if we define $L(-1)V=\span\{L(-1)v : v \in V\}$, we can consider the quotient space $V/L(-1)V$.  This has the structure of a commutative algebra under the operation $u \cdot v= \Res_x x^{-1} Y(u,x)v= u_{-1}v$.

%%%%%%%%%%%%%%%%%%%%%
\section{Modules}
%%%%%%%%%%%%%%%%%%%%%

There are a wide variety of definitions of modules for vertex operator algebras.  This variety stems from the amount of grading assumed for a given module and finite-dimensionality of the graded pieces (or lack thereof).  Some modules are ungraded and others admit a grading by $\N$, $\Q$, $\R$, or $\C$.   A $\N$-grading emphasizes lower-truncation, while the other gradings are given by the $L(0)$-eigenvalues.  With a grading in place we may impose a further restriction: the graded pieces must be finite-dimensional.  

Not only are there a variety of definitions for modules, the situation is further muddled by different names for the same objects (e.g.,``$\N$-graded weak '' and ``admissible'').  Other adjectives modifying ``module'' in the literature are: weak, strong, ordinary, lowest-weight, and generalized.  Because of the variety in language and structure, an explicit description of some of the different modules is warranted.  A natural starting point is modules for vertex algebras, which are naturally ungraded.  Every vertex operator algebra is a vertex algebra if one ignores the Virasoro vector and related axioms.

\begin{de} {\rm
For a vertex algebra $(V, Y, \1)$,  a {\it vertex algebra module} $(M, Y_M)$ is a vector space $M$ with a linear map

\begin{eqnarray}
Y_M:  &&V \rightarrow \End(M)[[x,x^{-1}]]\\
&&  v \mapsto Y_M(v,x)=\sum_{n \in \Z}v_n x^{-n-1}.
\end{eqnarray}

In addition $Y_M$ satisfies the following:

1) $v_nw=0$ for $n>>0$ where $v \in V$ and $w \in M$

2) $Y_M( \1,x)=Id_M$

3) For all $u,v \in V$,
\begin{eqnarray}
& &x_0^{-1}\d \left({x_1 - x_2 \over x_0}\right)Y_M(u,x_1)Y_M(v,x_2)-
x_0^{-1} \d \left({x_2- x_1 \over -x_0}\right)Y_M(v,x_2)Y_M(u,x_1) \nonumber \\
& &\ \ \ \ \ \ \ \ \ \ =x_2^{-1} \d \left({x_1- x_0 \over
x_2}\right)Y_M(Y(u,x_0)v,x_2).
\end{eqnarray}}
\end{de}

For a vertex operator algebra $V=(V,Y,\omega,\1)$, we can consider objects $(M, Y_M)$ as defined above for the vertex algebra structure of $V$. 

\begin{de} {\rm
A  {\it weak module} for a vertex operator algebra $V$ is a vertex algebra module for the vertex algebra structure of $V$.}
\end{de}

Weak modules for vertex operator algebras have additional structure that is a consequence of the vertex algebra module axioms.  They admit a representation of the Virasoro algebra and modules for a vertex operator also obey the $L(-1)$-derivation property.

\begin{prop}
Let $V=(V,Y, \omega, \1)$ be a vertex operator algebra and $M=(M, Y_M)$ a weak module for $V$.  \\
1) $Y_M(\om,x)=\sum_{n \in \Z}L_M(n) x^{-n-2}$ where \\
$$[L_M(m),L_M(n)]=(m-n)L_M(m+n)+\frac{m^3-m}{12}\d_{m+n,0}c$$
2) $Y_M(L(-1)v,x)=\frac{d}{dx}Y_M(v,x)$ for all $v \in V$
\end{prop}
 
Even with this additional structure known, weak modules of vertex operator algebras still lack suitable structure.  Some grading is necessary, and in particular a lower-truncated grading is desirable.  A lower-truncated grading guarantees the existence of ``lowest weight'' vectors.

\begin{de} {\rm
A weak module $M$ for a vertex operator algebra $V$ is called {\it $\N$-gradable} if it admits an $\N$-grading, $M=\bigoplus_{n \in \N} M(n)$, such that if $v \in V_r$ then $v_m M(n) \subseteq M(n+r-m-1).$}
\end{de}

The additional structure we have imposed on these modules is a lower-truncated grading, and we ensure that the grading is compatible with the vertex operator algebra action.  These modules are also called ``admissible'' in the literature.  The grading of these $\N$-gradable weak modules differs from the grading of vertex operator algebras in the following way.  The grading of vertex operator algebra is given by the eigenvalues of $L(0)$, while this is not necessarily true for for $\N$-gradable modules.  A third type of module is one where the grading is given by the $L(0)$-action. 

\begin{de} {\rm
A weak module $M$ for a  vertex operator algebra $V=(V,Y,\1,\omega)$ is a {\it $V$-module} if $M$ is $\C$-graded with $M=\bigoplus_{\l \in \C} M_{\l}$, and

1) $\dim(M_{\l})< \infty$,

2) $M_{\l+n=0}$ for fixed $\l$ and $n<<0$,

3) $L(0)w=\l w=\wt(w) w$, for $w \in M_{\l}$}.
\end{de}

The grading has been expanded to $\C$ to account for all possible $L(0)$-eigenvalues, and there is a lower truncation condition.  In addition, each graded piece must be finite-dimensional.  Such a finiteness condition is not imposed on $\N$-gradable weak modules.  One result of this finiteness condition and lower-truncation condition for $V$-modules is that $V$-modules are $\N$-gradable weak modules.  In practice, $\N$-gradable weak modules have enough structure to develop interesting theory.  We will see that for $C_2$-cofinite vertex operator algebras, weak modules are $\N$-gradable as well.  In his work on modularity, Zhu used what he called strong modules.  The definition of a strong module is the same as the definition of an ordinary module except that the axiom ``$\dim(M_{\l})< \infty$'' is omitted.

It is possible to extend the definition of $C_n$-cofiniteness to modules using $C_n(M)=\span\{\Res_x x^{-n} Y_M(u,x)w | u \in V, w\in M$\} for $n \geq 2 $.  Because there is no creation axiom for modules, it can be interesting to extend the idea of $C_1$-cofiniteness to modules in the naive way.

\begin{de} {\rm
For a vertex operator algebra $V=\bigoplus_{n \geq 0}V_n$ and a module $M$, let $$c_1(M)=\span\{\Res_x x^{-1} Y_M(u,x)w: u \in \bigoplus_{n>0}V_n, w \in M\}.$$  Then $M$ is called {\it $c_1$-cofinite} if $M/c_1(M)$ is finite-dimensional.}
\end{de}

This definition appears in the work of Nahm, who studied vertex operator algebras for which all irreducible $\N$-gradable weak modules are $c_1$-cofinite \cite{MR1305167}.  He called such algebras {\it quasirational}.  Quasirationality or $c_1$-cofiniteness of all irreducible modules is a important assumption in Huang's work on modular tensor categories and the Verlinde conjecture \cite{MR2029793} \cite{MR2151865} \cite{MR2140309}.  Huang's work also requires that the algebras be $C_2$-cofinite, which implies $c_1$-cofiniteness of the modules \cite{MR2052955}. 

%%%%%%%%%%%%%%%%%%%%%%%%%%%%
\section{Complete reducibility}
%%%%%%%%%%%%%%%%%%%%%%%%%%%%

One desirable property of vertex operator algebras that is featured in both mathematics and physics is complete reducibility of modules, the primary feature of ``rationality''.  The definition differs from author to author, with each rendition of ``rationality'' encompassing some minimum amount of ``goodness'' needed for the author's theory to work.  The ``goodness'' invariably includes some form of complete reducibility of modules and may also include some finiteness condition, i.e., finite number of irreducible modules, the graded pieces of irreducible modules are finite-dimensional, or even $C_2$-cofiniteness for some authors.  Calling a vertex operator algebra or conformal field theory ``rational'' endows that object with some physical importance, but the cost can sometimes be misinterpretation.  A common type of complete reducibility imposed on a vertex operator algebras is the following.\\

{\it Every $\N$-gradable weak module is the direct sum of irreducible $\N$-gradable weak modules.}\\

In mathematical literature, this property is sometimes called rationality, but certainly not consistently. A clearer naming would be complete reducibility of $\N$-gradable weak modules (in terms of irreducible $\N$-gradable weak modules).   I will use ``complete reducibility of $\N$-gradable  weak modules'' to convey this form of complete reducibility.  The assumption of this form of complete reducibility is necessary to prove many important results in vertex operator algebra theory.  As mentioned above,  the concept of rationality sometimes includes some finiteness assumptions.  Zhu's formulation of rationality included two additional conditions:  there exists a finite number of irreducible $\N$-gradable weak modules, and each graded piece of an irreducible $\N$-gradable weak module is finite-dimensional.  However Dong, Li, and Mason demonstrated that Zhu's additional conditions are consequences of complete reducibility of $\N$-gradable modules \cite{MR1628239}.  In other words, Zhu's seemingly stronger formulation of rationality is equivalent to complete reducibility of $\N$-gradable weak modules.   In fact, the Dong-Li-Mason results imply that complete reducibility of $\N$-gradable weak modules is equivalent to: every $\N$-gradable weak module is the direct sum of irreducible $V$-modules.  Some vertex operator algebras feature a stronger form of complete reducibility:\\

{\it Every weak module is the direct sum of irreducible $V$-modules.}\\

This property is called {\it regularity}, and examples of vertex operator algebras that satisfy this form of complete reducibility are the Moonshine module vertex operator algebra $V^{\natural}$, the Virasoro vertex operator algebras $L(c_{p,q},0)$, and vertex operator algebras associated to positive definite even lattices \cite{MR1488241}.  We will see that many more vertex operator algebras are regular.  Zhu conjectured that complete reducibility of $\N$-gradable modules implies $C_2$-cofiniteness.  This remains an important open question.  However for the stronger form of complete reducibility, Li proved that regular vertex operator algebras are $C_2$-cofinite \cite{MR1676852}.

%%%%%%%%%%%%%%%%%%%%%%%
\section{Spanning sets for algebras and modules}
%%%%%%%%%%%%%%%%%%%%%%%

One of the important consequences of $C_1$- or $C_2$-cofiniteness for a vertex operator algebra is that the algebra is finitely generated and has a Poincar\'{e}-Birkhoff-Witt-like spanning sets featuring desirable ordering restrictions.

\begin{prop} (cf. \cite{MR2023933})
For a subset $S$ of a vertex operator algebra $V=(V,Y, \om,\1)$, the subalgebra of $V$ generated by $S$ is
$$\<S\>=\span\{ u^{(1)}_{n_1} \cdots u^{(r)}_{n_r}\1 | r \in \N, u^{(1)}, \ldots, u^{(r)} \in S\cup\{\om\}, n_1, \ldots , n_r \in \Z\}.$$
\end{prop}

Different types of spanning sets feature different restrictions on the basic form, $u^{(1)}_{n_1} \cdots u^{(r)}_{n_r}\1$, of spanning set elements.  Some restrictions describe how often the index $n_i$ of or the weight of a mode $u^{(i)}_{n_i}$ can appear in a spanning set element, while other restrictions limit the $u^{(i)}$'s to certain subsets of $V$.   One way to think about the index restrictions on spanning set elements is in terms of a difference condition, similar to a difference condition on partitions.  A difference-$n$ condition on modes means that the indices of adjacent modes must differ by at least $n$. That is, for adjacent modes $u^{(i)}_{m_i}$ and $u^{(i+1)}_{m_{i+1}}$ in a spanning set element, $m_{i+1}-m_i \geq n$.

A natural question is: for a vertex operator algebra $V$, what sets $S$ generate $V$?  Certainly a minimal set $S$ is desirable, and this is what Karel and Li have explored.  

\begin{prop} (cf. \cite{MR1700507})
For a vertex operator algebra $V$, let $X$ be a set of homogeneous representatives of a spanning set for the quotient space $V/C_1(V)$. Then $V$ is spanned by the elements of the form
$$u^{(1)}_{n_1} \cdots u^{(r)}_{n_r}\1,$$
where $r \in \N$, $u^{(1)}, \ldots, u^{(r)} \in X$, $n_1, \ldots , n_r \in \Z$, and $\wt(u^1_{n_1}) \geq \cdots \geq \wt(u^r_{n_r}) >0$.
\end{prop}

In addition to showing that representatives of a basis for $V/C_1(V)$ generate $V$, Karel and Li also show that this set is a minimal generating  set of $V$.  So $C_1$-cofinite vertex operator algebras are finitely generated.  Karel and Li also prove an analogous spanning set for $\N$-gradable weak modules.

\begin{prop} (cf. \cite{MR1700507})
For a vertex operator algebra $V$ and an irreducible $\N$-gradable weak module $M=\bigoplus_{n \geq 0}M(n)$, let $X$ be a set of homogeneous representatives of a spanning set for the quotient space $V/C_1(V)$. Then $M$ is spanned by the elements of the form
$$u^{(1)}_{n_1} \cdots u^{(r)}_{n_r}w,$$
where $r \in \N$, $u^{(1)}, \ldots, u^{(r)} \in X$, $n_1, \ldots , n_r \in \Z$, $w \in M(0)$, and $\wt(u^1_{n_1}) \geq \cdots \geq \wt(u^r_{n_r}) >0$.
\end{prop}

For both the algebra and module spanning sets, the ordering restriction on the modes $u^{(i)}_{n_i}$ is in terms of the weight of the mode, and there is no restriction on how often an index $n_i$ of a mode $u^{(i)}_{n_i}$ can appear in a spanning set element.  However, one can prove an alternate version of the algebra spanning set, where the ordering restriction on the operators is in terms of the indices of modes, i.e., the $n_i$'s.

\begin{prop}
For a vertex operator algebra $V$, let $X$ be a set of homogeneous representatives of a spanning set for the quotient space $V/C_1(V)$. Then $V$ is spanned by the elements of the form
$$u^{(1)}_{n_1} \cdots u^{(r)}_{n_r}\1,$$
where $r \in \N$, $u^{(1)}, \ldots, u^{(r)} \in X$, $n_1, \ldots , n_r \in \Z$, and $n_1 \leq \cdots \leq {n_r}<0$. \end{prop}

The proof of this algebra spanning set is the same as the proof of the Karel-Li spanning result, since the mechanism for reordering the modes is the same.  This same  mechanism is just used to impose a different ordering.  It is possible to extend this spanning set to modules.    Since there is no creation axiom for modules, modes $u_n$ with $n \geq 0$ need to be limited in some way in the expression of spanning set elements. 

\begin{lem}
Given an $\N$-gradable weak module $M=\bigoplus_{n \geq 0} M(n)$ and $X$ a finite set of vectors in $V$, there exists $T \in \N$ such that $u_n w=0$ for all $n\geq T$, $u \in X$, and $w \in M(0)$.
\end{lem}

\pf We have $u_{\wt u +L} w = 0$ for all $v \in V$, $L \geq 0$, and $w \in M(0)$.  Let $T=\max_{u \in X}\{\wt u\}$. \qed

In particular, if $X$ is a set of representatives of a basis of $V/C_1(V)$ for a $C_1$-cofinite vertex operator algebra, such a $T$ exists.

\begin{prop}
For a $C_1$-cofinite vertex operator algebra $V$ and an irreducible $\N$-gradable weak module $M=\bigoplus_{n \geq 0}M(n)$, let $X$ be a set of homogeneous representatives of a spanning set for the quotient space $V/C_1(V)$. Then $M$ is spanned by the elements of the form
$$u^{(1)}_{n_1} \cdots u^{(r)}_{n_r}w,$$
where $r \in \N$, $u^{(1)}, \ldots, u^{(r)} \in X$, $n_1, \ldots , n_r \in \Z$, $w \in M(0)$, and $n_1 \leq \cdots \leq {n_r}<T$ (as above). \end{prop}

Gaberdiel and Neitzke developed another type of spanning set for a vertex operator algebra using a set of representatives of a basis of the quotient space $V/C_2(V)$.  Though this generating set is not minimal, it does have stronger ordering restrictions than the spanning set of Karel and Li.

\begin{prop} \label{gn} (cf. \cite{MR1990879}) 
For a vertex operator algebra $V$, let $X$ be a set of homogeneous representatives of a spanning set for the quotient space $V/C_2(V)$. Then $V$ is spanned by the elements of the form
$$u^{(1)}_{n_1} \cdots u^{(r)}_{n_r}\1,$$
where $r \in \N$, $u^{(1)}, \ldots, u^{(r)} \in X$, $n_1, \ldots , n_r \in \Z$, and $n_1 < \cdots < {n_r}<0$. 
\end{prop}

By enlarging the generating set, Gaberdiel and Neitzke were able to introduce a repetition restriction.  Each index of a mode can only appear once in the expression of a spanning set element, or in other words this is a no-repetion restriction on the indices of modes.  One corollary of Proposition \ref{gn} is that $C_2$-cofiniteness implies $C_n$-cofiniteness for $n \geq 2$.  The converse, mentioned above, is also true, yielding the following result.

\begin{cor}  (cf. \cite{MR1990879})
If a vertex operator algebra $V$ is $C_n$-cofinite for some $n \geq 2$ then $V$ is $C_n$-cofinite for all $n \geq 2$.
\end{cor}

$C_1$-cofiniteness of a vertex operator algebra is a strictly weaker condition since the vertex operator algebra constructed from a Heisenberg algebra is $C_1$-cofinite, but is not $C_2$-cofinite.

A more natural way to view the no-repetition restriction is in terms of a difference condition.  The Gaberdiel and Neitzke algebra spanning set obeys  a difference-one condition, and the reformulation of the  Karel and Li algebra spanning set obeys a difference-zero condition.  A natural extension of the Gaberdiel-Neitzke result would be a module spanning set satisfying a difference-one condition.  This next result is a partial solution to this difference-one module spanning set question.

\begin{prop} (cf. \cite{MR1927435})
\label{mythm} For a $C_2$-cofinite vertex operator algebra $V$ and an irreducible $\N$-gradable weak module $M=\bigoplus_{n \geq 0}M(n)$, let $X$ be a set of homogeneous representatives of a spanning set for the quotient space $V/C_2(V)$. Then $M$ is spanned by the elements of the form
$$u^{(1)}_{n_1} \cdots u^{(r)}_{n_r} w,$$
where $r \in \N$, $u^{(1)}, \ldots, u^{(r)} \in X$, $n_1, \ldots , n_r \in \Z$, $w \in M(0)$, and $u^1_{n_1}\leq \cdots \leq u^r_{n_r}<T$ (with $T$ as above) where $n_{j-1}<n_{j}$ if $n_j<0$ and $n_{j}=n_{j+1}$ for at most $Q$ indicies $j$ for $n_j \geq 0$, where $Q\in \N$, and $Q$ is fixed for $V$.
\end{prop}

In this module spanning set, the modes with negative indices obey a difference-one condition, but the non-negative modes do not.  However, the non-negative modes may repeat only  a globally finite number of times. This spanning set was useful in proving a number of results, yet it still is not a true difference-one condition module spanning set.  Miyamoto provides a further refinement obtaining a full difference-one module spanning set.

\begin{prop}  (cf. \cite{MR2046807})
For a $C_2$-cofinite vertex operator algebra $V$ and an irreducible $\N$-gradable weak module $M=\bigoplus_{n \geq 0}M(n)$, let $X$ be a set of homogeneous representatives of a spanning set for the quotient space $V/C_2(V)$. Then $M$ is spanned by the elements of the form
$$u^{(1)}_{n_1} \cdots u^{(r)}_{n_r}w,$$
where $r \in \N$, $u^{(1)}, \ldots, u^{(r)} \in X$, $n_1, \ldots , n_r \in \Z$, $w \in M(0)$, and $u^1_{n_1}< \cdots < u^r_{n_r}<T$ (as above).
\end{prop}

Viewed in terms of difference conditions, this means that $C_1$-cofiniteness implies a difference-zero condition on elements of a spanning set of a vertex operator algebra and its modules, and  $C_2$-cofiniteness implies a difference-one condition on elements of a spanning set of a vertex operator algebra and its modules.  

Orbifold theory and twisted modules are important aspects of the representation theory of vertex operator algebras.  A paper by Yamauchi \cite{MR2039213}
 addresses twisted modules, and the full statement of his difference-one spanning set theorem applies in this larger generality.

Again an underlying assumption in this exposition is that the vertex operator algebras are of CFT-type.  Miyamoto's result is true for vertex operator algebras that are not of CFT-type \cite{MR2046807}.  In particular, he assumed that $V=\bigoplus_{n\geq 0}V_n$, but $V_0$ not necessarily  one-dimensional.

%%%%%%%%%%%%%%%%%%%%%%%
\section{Finiteness Results}
%%%%%%%%%%%%%%%%%%%%%%%

%implications of C_2

As mentioned in previous sections, $C_2$-cofiniteness implies the finite-dimensionality of many quotient spaces of the algebra and implies the existence of a finite generating set for the algebra.  The assumption of $C_2$-finiteness has implications for the representation theory of vertex operator algebras beyond the difference-one module spanning set.

\begin{thm} \label{firstbig}
Let $V$ be a $C_2$-cofinite vertex operator algebra. Then:\\
1. $V$ has finite number of irreducible $V$-modules up to isomorphism. \cite{MR1700507}\\
2. Weak modules for $V$ are $\N$-gradable weak modules.\cite{MR2052955}\\
3. Irreducible $\N$-gradable weak modules  for $V$ are irreducible $V$-modules.\cite{MR1700507}\\
4. Irreducible weak modules for $V$ are  irreducible $V$-modules. \cite{MR2052955}\\
5. The associative algebra $A(V)$ is finite-dimensional.\\
\end{thm}

Practically, this means that under the assumption of $C_2$-cofiniteness, we do not need to be concerned about the myriad types of modules.  The weakest definition of modules is sufficient, as weak modules are gradable and lower truncated.  Further, any irreducible module has a grading given by the $L(0)$-action and each graded piece is finite-dimensional.  Some of these results were extended by Miyamoto, in his extended generality described above.

\begin{thm} \cite{MR2046807}
For $V$ a vertex operator algebra, the following are equivalent:\\
1. $V$ is $C_2$-cofinite.\\
2. Every weak module is a direct sum of generalized eigenspaces of $L(0)$.\\
3. Every weak module is an $\N$-gradable weak module $M = \bigoplus_{n \geq 0} M(n)$ such that $M(n)$ is a direct sum of generalized eigenspaces of $L(0)$.\\
4. $V$ is finitely generated and every weak module is an $\N$-gradable weak module.
 \end{thm}

In light of this, we see that $C_2$-cofiniteness is equivalent to all modules having suitable properties for an interesting representation theory, with the lone exception of complete reducibility.  However, the assumption of $C_2$-cofiniteness unifies notions of complete reducibility.

\begin{thm}\cite{MR2052955} \cite{MR1676852}
For a $C_2$-cofinite vertex operator algebra $V$, the following are equivalent:\\
1. Every weak module for $V$ is the direct sum of irreducible $V$-modules.\\
2. Every $\N$-gradable weak module is the direct sum of irreducible $\N$-gradable weak modules.
\end{thm}

In particular, this means that all known vertex operator algebras with complete reducibility of $\N$-gradable weak modules are regular.  Theorem \ref{firstbig} should be compared with the following theorem for vertex operator algebras with complete reducibility of $\N$-gradable weak modules.

\begin{thm}
Let $V$ be a vertex operator algebra for which every $\N$-gradable weak module is the direct sum of irreducible $V$-modules. Then:\\
1. $V$ has a finite number of irreducible $V$ modules up to isomorphism. \cite{MR1615132}\\
2. Irreducible $\N$-gradable weak modules  for $V$ are irreducible $V$-modules. \cite{MR1615132}\\
3. The associative algebra $A(V)$ is semisimple and finite-dimensional.\cite{MR1637252}\\
\end{thm}

This is compelling evidence that complete reducibility of $\N$-gradable weak modules and $C_2$-cofiniteness are somehow related.  Zhu conjectured that complete reducibility of $\N$-gradable weak modules implies $C_2$-cofiniteness \cite{MR1317233}. The converse of this conjecture has been disproved. Building on the work of Kausch and Gaberdiel \cite{MR1411388}, Abe and Carqueville-Flohr construct examples of $C_2$-cofinite vertex operators for which there exist $\N$-gradable weak modules that are not completely reducible.  Specifically, Abe constructs a family of $C_2$-cofinite vertex operator algebras with central charge $-2d$ for $d \in \Z_+$ with reducible indecomposible modules \cite{math.QA/0503472}. Carqueville and Flohr prove that the vertex operator algebras constructed from the triplet algebras $c_{p,1}$ are $C_2$-cofinite and also have reducible indecomposible modules \cite{carqueville-2006-39}.  However Zhu's conjecture as to whether complete reducibility of $\N$-gradable weak modules implies $C_2$-cofiniteness remains open.

%%%%%%%%%%%%%%%%%%%%%%%%%%%%%
%bibliography
%%%%%%%%%%%%%%%%%%%%%%%%%%%%%

%\bibliographystyle{unsrt}   % this means that the order of references
			    % is dtermined by the order in which the
			    % \cite and \nocite commands appear
%\bibliography{gbuhl_edinburgh}  % list here all the bibliographies that

\begin{thebibliography}{99}

\bibitem{math.QA/0503472}
T.~Abe.
\newblock {A $Z_2$-orbifold model of the symplectic fermionic vertex operator
  superalgebra}.
\newblock arXiv:math.QA/0503472.

\bibitem{MR2052955}
T.~Abe, G.~Buhl, and C.~Dong.
\newblock Rationality, regularity, and {$C\sb 2$}-cofiniteness.
\newblock {\em Trans. Amer. Math. Soc.}, 356(8):3391--3402, 2004.

\bibitem{MR1172696}
R.~E. Borcherds.
\newblock Monstrous moonshine and monstrous {L}ie superalgebras.
\newblock {\em Invent. Math.}, 109(2):405--444, 1992.

\bibitem{MR1927435}
G.~Buhl.
\newblock A spanning set for {VOA} modules.
\newblock {\em J. Algebra}, 254(1):125--151, 2002.

\bibitem{carqueville-2006-39}
N.~Carqueville and M.~Flohr.
\newblock Nonmeromorphic operator product expansion and $c_2$-cofiniteness for
  a family of w-algebras.
\newblock {\em J.PHYS.A}, 39:951, 2006.

\bibitem{MR1488241}
C.~Dong, H.~Li, and G.~Mason.
\newblock Regularity of rational vertex operator algebras.
\newblock {\em Adv. Math.}, 132(1):148--166, 1997.

\bibitem{MR1615132}
C.~Dong, H.~Li, and G.~Mason.
\newblock Twisted representations of vertex operator algebras.
\newblock {\em Math. Ann.}, 310(3):571--600, 1998.

\bibitem{MR1628239}
C.~Dong, H.~Li, and G.~Mason.
\newblock Twisted representations of vertex operator algebras and associative
  algebras.
\newblock {\em Internat. Math. Res. Notices}, (8):389--397, 1998.

\bibitem{MR1637252}
C.~Dong, H.~Li, and G.~Mason.
\newblock Vertex operator algebras and associative algebras.
\newblock {\em J. Algebra}, 206(1):67--96, 1998.

\bibitem{MR996026}
I.~Frenkel, J.~Lepowsky, and A.~Meurman.
\newblock {\em Vertex operator algebras and the {M}onster}, volume 134 of {\em
  Pure and Applied Mathematics}.
\newblock Academic Press Inc., Boston, MA, 1988.

\bibitem{MR1411388}
M.~R. Gaberdiel and H.~G. Kausch.
\newblock A rational logarithmic conformal field theory.
\newblock {\em Phys. Lett. B}, 386(1-4):131--137, 1996.

\bibitem{MR1990879}
M.~R. Gaberdiel and A.~Neitzke.
\newblock Rationality, quasirationality and finite {$W$}-algebras.
\newblock {\em Comm. Math. Phys.}, 238(1-2):305--331, 2003.

\bibitem{MR2029793}
Y.-Z. Huang.
\newblock Riemann surfaces with boundaries and the theory of vertex operator
  algebras.
\newblock In {\em Vertex operator algebras in mathematics and physics (Toronto,
  ON, 2000)}, volume~39 of {\em Fields Inst. Commun.}, pages 109--125. Amer.
  Math. Soc., Providence, RI, 2003.

\bibitem{MR2151865}
Y.-Z. Huang.
\newblock Differential equations and intertwining operators.
\newblock {\em Commun. Contemp. Math.}, 7(3):375--400, 2005.

\bibitem{MR2140309}
Y.-Z. Huang.
\newblock Vertex operator algebras, the {V}erlinde conjecture, and modular
  tensor categories.
\newblock {\em Proc. Natl. Acad. Sci. USA}, 102(15):5352--5356 (electronic),
  2005.

\bibitem{math.QA/0311235}
Y.-Z. Huang, J.~Lepowsky, and L.~Zhang.
\newblock {A logarithmic generalization of tensor product theory for modules
  for a vertex operator algebra}.
\newblock {\em Internat. J. Math., to appear}, arXiv:math.QA/0311235.

\bibitem{MR1700507}
M.~Karel and H.~Li.
\newblock Certain generating subspaces for vertex operator algebras.
\newblock {\em J. Algebra}, 217(2):393--421, 1999.

\bibitem{MR2023933}
J.~Lepowsky and H.~Li.
\newblock {\em Introduction to vertex operator algebras and their
  representations}, volume 227 of {\em Progress in Mathematics}.
\newblock Birkh\"auser Boston Inc., Boston, MA, 2004.

\bibitem{MR1676852}
H.~Li.
\newblock Some finiteness properties of regular vertex operator algebras.
\newblock {\em J. Algebra}, 212(2):495--514, 1999.

\bibitem{MR2046807}
M.~Miyamoto.
\newblock Modular invariance of vertex operator algebras satisfying {$C\sb
  2$}-cofiniteness.
\newblock {\em Duke Math. J.}, 122(1):51--91, 2004.

\bibitem{MR1305167}
W.~Nahm.
\newblock Quasi-rational fusion products.
\newblock {\em Internat. J. Modern Phys. B}, 8(25-26):3693--3702, 1994.

\bibitem{MR2039213}
H.~Yamauchi.
\newblock Modularity on vertex operator algebras arising from semisimple
  primary vectors.
\newblock {\em Internat. J. Math.}, 15(1):87--109, 2004.

\bibitem{MR1317233}
Y.~Zhu.
\newblock Modular invariance of characters of vertex operator algebras.
\newblock {\em J. Amer. Math. Soc.}, 9(1):237--302, 1996.



\end{thebibliography}
			     % you need. 

\end{document}